\documentclass{amsart}
\usepackage{amsthm}
\usepackage{amsmath}
\usepackage{amsfonts}
\usepackage{amssymb}
\usepackage{mathrsfs}
\usepackage{hyperref}
\usepackage[numbers]{natbib}
\usepackage{fancyhdr}
\theoremstyle{plain}
\newtheorem{lem}{Lemma}[section]
\newtheorem{prop}[lem]{Proposition}
\newtheorem{thm}[lem]{Theorem}
\newtheorem{cor}[lem]{Corollary}

\theoremstyle{definition}
\newtheorem{defn}[lem]{Definition}
\newtheorem{rem}[lem]{Remark}
\newcommand{\R}{\mathbb R}
\newcommand{\N}{\mathbb N}
\newcommand{\Z}{\mathbb Z}
\newcommand{\Diff}{\mbox{\rm Diff}}
\newcommand{\Vect}{\mbox{\rm Vect}}
\newcommand{\C}{\mbox{\rm C}}
\newcommand{\id}{\mbox{\rm id}}
\newcommand{\dx}{\,\text{\rm d}x}
\renewcommand{\d}{\,\text{\rm d}}
\newcommand{\D}{\,\text{\rm D}}

\newcommand{\ad}{\mbox{\rm ad}}
\renewcommand{\S}{\mathbb S}
\renewcommand{\phi}{\varphi}
\newcommand{\muDP}{\mu\text{\rm DP}}
\newcommand{\norm}[1]{\left|\!\left|#1\right|\!\right|}
\newcommand{\eps}{\varepsilon}
\newcommand{\ska}[2]{\left\langle #1,#2\right\rangle}
\newcommand{\set}[2]{\left\{#1;\;#2\right\}}
\newcommand{\bea}{\begin{eqnarray}}
\newcommand{\eea}{\end{eqnarray}}
\begin{document}
\title{Geometric aspects of the periodic $\mu$DP equation}
\author{Joachim Escher, Martin Kohlmann, Boris Kolev}
\address{Institute for Applied Mathematics, University of Hannover, D-30167 Hannover, Germany}
\email{escher@ifam.uni-hannover.de, kohlmann@ifam.uni-hannover.de}
\address{CMI, 39 rue F. Joliot-Curie, 13453 Marseille Cedex 13,
France}
\email{kolev@cmi.univ-mrs.fr}
\begin{abstract} We consider the periodic $\muDP$ equation (a modified version
of the Degasperis-Procesi equation) as the geodesic flow of a
right-invariant affine connection $\nabla$ on the Fr\'echet Lie
group $\Diff^{\infty}(\S^1)$ of all smooth and
orientation-preserving diffeomorphisms of the circle $\S^1=\R/\Z$.
On the Lie algebra $\C^{\infty}(\S^1)$ of $\Diff^{\infty}(\S^1)$,
this connection is canonically given by the sum of the Lie bracket
and a bilinear operator. For smooth initial data, we show the
short time existence of a smooth solution of $\muDP$ which depends
smoothly on time and on the initial data. Furthermore, we prove
that the exponential map defined by $\nabla$ is a smooth local
diffeomorphism of a neighbourhood of zero in $\C^{\infty}(\S^1)$
onto a neighbourhood of the unit element in
$\Diff^{\infty}(\S^1)$. Our results follow from a general approach on non-metric Euler equations on Lie groups, a Banach space
approximation of the Fr\'echet space $\C^{\infty}(\S^1)$, and a
sharp spatial regularity result for the geodesic flow.
\end{abstract}
\maketitle
\begin{center}
{\em{Dedicated to Herbert Amann on the occasion of his 70th birthday}}\\[0.9cm]
\end{center}
\section{Introduction}
In recent years, several nonlinear equations arising as
approximations to the governing model equations for water waves
attracted a considerable amount of attention in the fluid dynamics
research community (cf.\ \cite{Johnson}). The Korteweg-de Vries
(KdV) equation is a well-known model for wave-motion on shallow
water with small amplitudes and a flat bottom. This equation is
completely integrable, allows for a Lax pair formulation and the
corresponding Cauchy problem was the subject of many studies.
However, it was observed in \cite{Bourgain} that solutions of the KdV equation do not
break as physical water waves do: The flow is globally well posed for square integrable initial data (see also
\cite{Kappeler,KenigPonceVega} for further results). The
Camassa-Holm (CH) equation
$$u_t+3uu_{x}=2u_{x}u_{xx}+uu_{xxx}+u_{txx}$$
was introduced to model the shallow-water medium-amplitude regime
(see \cite{CamassaHolm}). Closely related to the CH equation is the Degasperis-Procesi (DP) equation
$$u_t+4uu_{x}=3u_{x}u_{xx}+uu_{xxx}+u_{txx}$$
was discovered in a search for integrable equations similar to the
CH equation (see \cite{DegasperisProcesi}). Both equations are
higher order approximations in a small amplitude expansion of the
incompressible Euler equations for the unidirectional motion of
waves at a free surface under the influence of gravity (cf.\
\cite{ConstantinLannes}). They have a bi-Hamiltonian structure,
are completely integrable and allow for wave breaking and peaked
solitons, \cite{CE98b,DegasperisHolmHone,Escher08,Ivanov}. The
Cauchy problem for the periodic CH equation in spaces of classical
solutions has been studied extensively (see, e.g.,
\cite{ConstantinEscher2, Misiolek}); in \cite{ConstantinEscher1}
and \cite{DeLellis} the authors explain that this equation is also
well posed in spaces which include peakons, showing in this way
that peakons are indeed meaningful solutions of CH. Well-posedness
for the periodic DP equation and various features of solutions of
the DP on the circle are discussed in \cite{ELY}. Both, the CH
equation and the DP equation, are embedded into the family of
$b$-equations
\bea\label{beqn}m_t=-(m_xu+bmu_x),\quad m:=u-u_{xx},\eea
where $u(t,x)$ is a function of a spatial variable $x\in\S^1$ and
a temporal variable $t\in\mathbb{R}$. Note that the family (\ref{beqn}) can be
derived as the family of asymptotically equivalent shallow water
wave equations that emerges at quadratic order accuracy for any
$b\neq -1$ by an appropriate Kodama transformation,
\cite{DullinGottwaldHolm}. For $b=2$, we recover the CH equation
and for $b=3$, we get the DP equation. Note that the $b$-equation
is integrable only if $b=2$ or $b=3$. For further results and references we refer
to \cite{EscherYin}.\\\indent
Since the pioneering works \cite{Arnold,EbinMarsden}, geometric
interpretations of evolution equations led to several interesting
results in the applied analysis literature. A detailed discussion
of the CH equation in this framework was given by
\cite{Kouranbaeva}. The geometrical aspects of some metric Euler
equations are explained in
\cite{ConstantinKolev,ConstantinKolev2,Kolev,Michor}. Studying the
$b$-equations as a geodesic flow on the diffeomorphism group
$\Diff^{\infty}(\S^1)$, it was shown recently in \cite{EscherKolev} that for smooth
initial data $u_0(x)=u(0,x)$, there is a unique short-time
solution $u(t,x)$ of (\ref{beqn}), depending smoothly on
$(t,u_0)$. The crucial idea is to define an affine (not
necessarily Riemannian) connection $\nabla$ on
$\Diff^{\infty}(\S^1)$ given at the identity by the sum of
the Lie bracket and a bilinear symmetric operator $B$
so that $B(u,u)=-u_t$. Most importantly, this approach also works
for $b$-equations of non-metric type and it motivates the study of
geometric quantities like curvature or an exponential map for the
family (\ref{beqn}). In particular, the authors of \cite{EscherKolev} proved that
the exponential map for $\nabla$ is a smooth local diffeomorphism
near zero in $\C^{\infty}(\S^1)$. Recently it has been shown in \cite{EscherSeiler}
that the $b$-equation can be realized as a metric
Euler equation \emph{only} if $b=2$. In all other cases $b\ne 2$ there is no
Riemannian metric on $\Diff^{\infty}(\S^1)$ such that the
corresponding geodesic flow is re-expressed by the $b$-equation.
Geometric aspects of some novel nonlinear PDEs related to CH and
DP are discussed in \cite{LenellsMisiolekTiglay}.\\\indent
In this paper, we study the $\muDP$ equation
\bea\label{muDP}\mu(u_t)-u_{txx}+3\mu(u)u_x-3u_xu_{xx}-uu_{xxx}=0,\eea
where $\mu$ denotes the projection $\mu(u)=\int_0^1u\dx$ and
$u(t,x)$ is a spatially periodic real-valued function of a time
variable $t\in\R$ and a space variable $x\in\S^1$. The $\mu$DP equation
belongs to the family of $\mu$-$b$-equations which follows from
(\ref{beqn}) by replacing $m=\mu(u)-u_{xx}$. The study of
$\mu$-variants of (\ref{beqn}) is motivated by the following key
observation: Letting $m=-\partial_x^2u$, equation (\ref{beqn}) for
$b=2$ becomes the Hunter-Saxton (HS) equation
$$2u_xu_{xx}+uu_{xxx}+u_{txx}=0,$$
which possesses various interesting geometric properties, see,
e.g., \cite{Len1,Len2}, whereas the choice $m=(1-\partial_x^2)u$
leads to the CH equation as explained above. In the search for
integrable equations that are given by a perturbation of
$-\partial_x^2$, the $\mu$-$b$-equation has been introduced and it
could be shown that it behaves quite similarly to the
$b$-equation; see \cite{LenellsMisiolekTiglay} where the authors
discuss local and global well-posedness as well as finite time
blow-up and peakons. Our study of the $\muDP$ is inspired by the
results in \cite{EscherKolev}. In fact using the approach of \cite{EscherKolev}
we shall conceptualise a geometric picture of the $\muDP$ equation.\\
\indent Our study is mostly performed in the $\C^\infty$-category.
Elements of $\C^{\infty}(\S^1)$ are sometimes also called smooth
for brevity.

We will reformulate the $\muDP$
equation in terms of a geodesic flow on $\Diff^{\infty}(\S^1)$ to
obtain the following main result: Given a smooth initial data $u_0(x)$,
for which $\norm{u_0}_{\mbox{\rm\scriptsize C}^3(\S^1)}$ is small,
there is a unique smooth solution $u(t,x)$ of
(\ref{muDP}) which depends smoothly on $(t,u_0)$. More precisely, we have
\thm\label{cor_main_Escher_Kolev_09} There exists an open interval
$J$ centered at zero and $\delta>0$ such that for each
$u_0\in\C^{\infty}(\S^1)$ with $\norm{u_0}_{\mbox{\rm\scriptsize
C}^3(\S^1)}<\delta$, there exists a unique solution
$u\in\C^{\infty}(J,\C^{\infty}(\S^1))$ of the $\muDP$ equation
such that $u(0)=u_0$. Moreover, the solution $u$ depends smoothly
on $(t,u_0)\in J\times\C^{\infty}(\S^1)$.\endcor\rm
It is known that the Riemannian exponential mapping on general
Fr\'echet manifolds fails to be a smooth local diffeomorphism from
the tangent space back to the manifold, cf. \cite{ConstantinKolev}.
Therefore the following result is quite remarkable.

\thm\label{thm_exp} The exponential map $\exp$ at the unity
element for the $\muDP$ equation on $\Diff^{\infty}(\S^1)$ is a
smooth local diffeomorphism from a neighbourhood of zero in
$\C^{\infty}(\S^1)$ onto a neighbourhood of $\id$ in
$\Diff^{\infty}(\S^1)$.\endthm\rm
Our paper is organized as follows: In Section 2, we rewrite
(\ref{muDP}) in terms of a local flow $\phi\in\Diff^{n}(\S^1)$,
$n\geq 3$, and explain the geometric setting. The resulting
equation is an ordinary differential equation and in Section 3, we
apply the Theorem of Picard and Lindel\"of to obtain a solution of class
$\C^n(\S^1)$ with smooth dependence on $t$ and $u_0(x)$. In
addition, we show that this solution in
$\Diff^n(\S^1)\times\C^n(\S^1)$ does neither lose nor gain spatial
regularity as $t$ varies through the associated interval of existence. We
then approximate the Fr\'echet Lie group $\Diff^{\infty}(\S^1)$ by
the topological groups $\Diff^{n}(\S^1)$ and the Fr\'echet space
$\C^{\infty}(\S^1)$ by the Banach spaces $\C^n(\S^1)$ to obtain an
analogous existence result for the geodesic equation on
$\Diff^{\infty}(\S^1)$. Finally, in Section 4, we make again use
of a Banach space approximation to prove that the exponential map
for the $\muDP$ is a smooth local diffeomorphism nero
zero as a map $\C^{\infty}(\S^1)\to\Diff^{\infty}(\S^1)$.\\[.5cm]\indent
\textbf{Acknowledgements.} The authors would like to thank Mats
Ehrnstr\"om and Jonatan Lenells for several fruitful discussions
as well as the referee for helpful remarks.
%
%
%
%
\section{Geometric reformulation of the $\mu$DP equation}
We write $\Diff^{\infty}(\S^1)$ for the smooth
orientation-preserving diffeomorphisms of the unit circle
$\S^1=\R/\Z$ and $\Vect^{\infty}(\S^1)$ for the space of smooth
vector fields on $\S^1$. Clearly, $\Diff^{\infty}(\S^1)$ is a Lie
group and it is easy to see that its Lie algebra is
$\Vect^{\infty}(\S^1)$: If $t\mapsto\phi(t)$ is a smooth path in
$\Diff^{\infty}(\S^1)$ with $\phi(0)=\id$, then $\phi_t(0,x)\in
T_x\S^1$ for all $x\in\S^1$ and thus the Lie algebra element
$\phi_t(0,\cdot)$ is a smooth vector field on $\S^1$. Furthermore,
since $T\S^1\simeq\S^1\times\R$ is trivial, we can identify the
Lie algebra $\Vect^{\infty}(\S^1)$ with $\C^{\infty}(\S^1)$. Note
that $[u,v]=u_xv-v_xu$ is the corresponding Lie bracket. In the
following, we will also use that $\Diff^{\infty}(\S^1)$ has a
smooth manifold structure modelled over the Fr\'echet
space $\C^{\infty}(\S^1)$. In particular, $\Diff^{\infty}(\S^1)$
is a Fr\'echet Lie group and thus it is parallelizable, i.e.,
$T\Diff^{\infty}(\S^1)\simeq\Diff^{\infty}(\S^1)\times\C^{\infty}(\S^1)$.
Let $\Diff^n(\S^1)$ denote the group of orientation-preserving
diffeomorphisms of $\S^1$ which are of class $\C^n(\S^1)$.
Similarly, $\Diff^n(\S^1)$ has a smooth manifold structure
modelled over the Banach space $\C^n(\S^1)$. Note that
$\Diff^n(\S^1)$ is only a topological group but not a Banach Lie
group, since the composition and inversion maps are continuous but
not smooth. Furthermore, the trivialization
$T\Diff^n(\S^1)\simeq\Diff^n(\S^1)\times\C^n(\S^1)$ is only
topological and not smooth.\\\indent
In this section, we write (\ref{muDP}) as an ordinary differential
equation on the tangent bundle $\Diff^n(\S^1)\times\C^n(\S^1)$, where $n\geq 3$. In a
first step, we rewrite (\ref{muDP}) using the operator
$A:=\mu-\partial_x^2$. Here $\mu$ denotes the linear map given by
$f\mapsto\int_0^1f(t,x)\dx$ for any function $f(t,x)$ depending on
time $t$ and space $x\in\S^1$. Observe that
$\mu(\partial_x^kf)=0$, $k\geq 1$, if $f$ and its derivatives are
continuous functions on $\S^1$. Furthermore, $\mu(f)$ is still
depending on the time variable $t$. The following lemma
establishes the invertibility of $A$ as an operator acting on
$\C^n(\S^1)$ for $n\geq 2$.
\lem\label{lem_inv_of_A} Given $n\geq 2$, the operator $A=\mu-\partial_x^2$ maps
$\C^n(\S^1)$ isomorphically onto $\C^{n-2}(\S^1)$. The inverse is given by
\bea(A^{-1}f)(x)&=&\left(\frac{1}{2}x^2-\frac{1}{2}x+
\frac{13}{12}\right)\int_0^1f(a)\d a+\left(x-\frac{1}{2}\right)
\int_0^1\int_0^af(b)\d b\d a\nonumber\\
&&-\int_0^x\int_0^af(b)\d b\d a+\int_0^1\int_0^a\int_0^bf(c)\d c\d
b\d a.\nonumber\eea
\endlem\rm
\proof Clearly, $\mu(A^{-1}f)=\mu(f)$ and
$(A^{-1}f)_{xx}=\mu(f)-f$ so that $A(A^{-1}f)=f$. To verify that
$A$ is surjective, we observe that
$\partial_x^k(A^{-1}f)(0)=\partial_x^k(A^{-1}f)(1)$ for all
$k\in\{0,\dots,n\}$. To see that $A$ is injective, assume that
$Au=0$ for $u\in\C^n(\S^1)$ and $n\geq 2$. Then there are
constants $c,d\in\R$ such that $u=\frac{1}{2}\mu(u)x^2+cx+d$. By
periodicity we first conclude that $c=0$ and $\mu(u)=0$. Hence $d$
has to vanish as well.
\endproof
\lem\label{lem_muDP_rewritten} Assume that
$u\in\C((-T,T),\C^n(\S^1))\cap \C^1((-T,T),\C^{n-1}(\S^1))$ is a
solution of \text{\rm(\ref{muDP})} for some $n\geq 3$ with $T>0$.
Then the $\muDP$ equation can be written as
\bea\label{muDP_rewritten}
u_t=-A^{-1}(u(Au)_x+3(Au)u_x).\eea
\endlem
\proof Writing (\ref{muDP}) in the form
$$\mu(u_t)-u_{txx}=uu_{xxx}-3u_x(\mu(u)-u_{xx}),$$
we see that it is equivalent to
$$Au_t=-u(Au)_x-3(Au)u_x.$$
Thus $u$ is a solution of (\ref{muDP}) if and only if
(\ref{muDP_rewritten}) holds true.\endproof
As explained in \cite{Lang,LangMfds}, the vector field $u(t,x)$
admits a unique local flow $\phi$ of class $\C^n(\S^1)$, i.e.,
$$\phi_t(t,x)=u(t,\phi(t,x)),\quad\varphi(0,x)=x$$
for all $x\in\S^1$ and all $t$ in some open interval $J\subset\R$.
We will use the short-hand notation $\phi_t=u\circ\phi$ for
$\phi_t(t,x)=u(t,\phi(t,x))$; i.e., $\circ$ denotes the
composition with respect to the spatial variable. Particularly, we have that
$u=\phi_t\circ\phi^{-1}$. Moreover, given
$(\phi,\xi)\in\C^1(J,\Diff^n(\S^1)\times\C^n(\S^1))$,
then $\phi^{-1}(t)$ is a $\C^n(\S^1)$-diffeomorphism for all $t\in
J$ and $\xi\circ\phi^{-1}\in\C^1(J,\C^n(\S^1))$.\\\indent
In this paper, we are mainly interested in smooth diffeomorphisms
on $\S^1$. For the reader's convenience we briefly recall the basic
geometric setting.
Let us consider the Fr\'echet manifold $\Diff^{\infty}(\S^1)$ and
a continuous non-degenerate inner product $\ska{\cdot}{\cdot}$ on
$\C^{\infty}(\S^1)$, i.e., $u\mapsto\ska{u}{u}$ is continuous (and
hence smooth) and $\ska{u}{v}=0$ for all $v\in\C^{\infty}(\S^1)$
forces $u=0$. To define a weak right-invariant Riemannian metric
on $\Diff^{\infty}(\S^1)$, we extend the inner product
$\ska{\cdot}{\cdot}$ to any tangent space by right-translations,
i.e., for all $g\in\Diff^{\infty}(\S^1)$ and all $u,v\in
T_g\Diff^{\infty}(\S^1)$, we set
$$\ska{u}{v}_{g}=\ska{(R_{g^{-1}})_{*}u}{(R_{g^{-1}})_{*}v}_e,$$
where $e$ denotes the identity. Observe that any open set in the
topology induced by this inner product is open in the Fr\'echet
space topology of $\C^{\infty}(\S^1)$ but the converse is not
true. We therefore call $\ska{\cdot}{\cdot}$ a \emph{weak} Riemannian
metric on $\Diff^{\infty}(\S^1)$, cf. \cite{ConstantinKolev}. We
next define a bilinear operator
$B:\Vect^{\infty}(\S^1)\times\Vect^{\infty}(\S^1)\to\Vect^{\infty}(\S^1)$
by
$$B(u,v)=\frac{1}{2}((\ad_u)^*(v)+(\ad_v)^*(u)),$$
where $(\ad_u)^*$ is the adjoint (with respect to
$\ska{\cdot}{\cdot}$) of the natural action of the Lie algebra on
itself given by $\ad_u:v\mapsto[u,v]$. Observe that $B$ defines a
right-invariant affine connection $\nabla$ on
$\Diff^{\infty}(\S^1)$ by
\bea\label{connection}\nabla_{\xi_u}\xi_v=\frac{1}{2}[\xi_u,\xi_v]+B(\xi_u,\xi_v),\eea
where $\xi_u$ and $\xi_v$ are the right-invariant vector fields on
$\Diff^{\infty}(\S^1)$ with values $u,v$ at the identity. It can be shown
that a smooth curve $t\mapsto g(t)$ in $\Diff^{\infty}(\S^1)$ is a
geodesic if and only if $u=(R_{g^{-1}})_*\dot g$ solves the
\emph{Euler equation}
\bea\label{Euler} u_t=-B(u,u);\eea
here, $u$ is the \emph{Eulerian velocity} (cf.\
\cite{ArnoldMech}). Hence the Euler equation (\ref{Euler})
corresponds to the geodesic flow of the affine connection $\nabla$
on the diffeomorphism group $\Diff^{\infty}(\S^1)$. Paradigmatic
examples are the following: In \cite{ConstantinKolev}, the authors
show that the Euler equation for the right-invariant $L^2$-metric
on $\Diff^{\infty}(\S^1)$ is given by the inviscid Burgers
equation. Equipping on the other hand $\C^{\infty}(\S^1)$ with the
$H^1$-metric, one obtains the Camassa-Holm equation. Similar
correspondences for the general $H^k$-metrics are explained in
\cite{ConstantinKolev2}.\\
\indent Conversely, starting with an equation of type
$u_t=-B(u,u)$ with a bilinear operator $B$, one associates an
affine connection $\nabla$ on $\Diff^{\infty}(\S^1)$ by formula
(\ref{connection}). It is however by no means clear that this
connection corresponds to a Riemannian structure on
$\Diff^{\infty}(\S^1)$. It is worthwhile to mention that the
connection $\nabla$ corresponding to the family of $b$-equations
is compatible with some metric only for $b=2$: In
\cite{EscherSeiler} the authors explain that for any $b\neq 2$,
the $b$-equation (\ref{beqn}) cannot be realized as an Euler
equation on $\Diff^{\infty}(\S^1)$ for any regular inertia
operator. This motivates the notion of \emph{non-metric Euler
equations}. An analogous result holds true for the
$\mu$-$b$-equations from which we conclude that the $\muDP$
equation belongs to the class of non-metric Euler equations.
Although we have no metric for the $\muDP$ equation, we will
obtain some geometric information by using the connection
$\nabla$, defined in the following way.\\\indent
Let $X(t)=(\phi(t),\xi(t))$ be a vector field along the curve
$\phi(t)\in\Diff^{\infty}(\S^1)$. Furthermore let
$$B(v,w):=\frac{1}{2}A^{-1}(v(Aw)_x+w(Av)_x+3(Av)w_x+3(Aw)v_x).$$
Lemma \ref{lem_muDP_rewritten} shows that
$$B(u,u)=A^{-1}(u(Au)_x+3(Au)u_x)=-u_t,$$
if $u$ is a solution to the $\muDP$ equation. Next,
the covariant derivative of $X(t)$ in the present case is defined
as
$$\frac{\!\D X}{\!\D t}(t)=\left(\phi(t),\xi_t+
\frac{1}{2}[u(t),\xi(t)]+B(u(t),\xi(t))\right),$$
where $u=\phi_t\circ\phi^{-1}$. We see that $u$ is a solution of
the $\muDP$ if and only if its local flow $\phi$ is a geodesic for
the connection $\nabla$ defined by $B$ via
(\ref{connection}).\\\indent
Although we are mainly interested in the smooth category, we
will first discuss flows $\phi(t)$ on $\Diff^{n}(\S^1)$ for
technical purposes. Regarding $\Diff^n(\S^1)$ as a smooth Banach
manifold modelled over $\C^n(\S^1)$, the following result has to be
understood locally, i.e., in any local chart of $\Diff^n(\S^1)$.
\prop\label{prop_muDP_geo} Given $n\geq 3$, the function
$u\in\C(J,\C^n(\S^1))\cap\C^1(J,\C^{n-1}(\S^1))$
is a solution of \text{\rm(\ref{muDP})} if and only if
$(\phi,\xi)\in\C^1(J,\Diff^n(\S^1)\times\C^n(\S^1))$ is a solution
of
\bea\label{muDPgeo}\left\{\begin{array}{ccc}
  \phi_t & = & \xi, \\
  \xi_t & = & -P_{\phi}(\xi),\end{array}\right.\eea
where $P_{\phi}:=R_{\phi}\circ P\circ R_{\phi^{-1}}$ and
$P(f):=3A^{-1}(f_xf_{xx}+(Af)f_x)$.
\endprop
\proof The function $u$ and the corresponding flow
$\phi\in\Diff^n(\S^1)$ satisfy the relation $\phi_t=u\circ\phi$.
Setting $\phi_t=\xi$, the chain rule implies that
$$\xi_t=(u_t+uu_x)\circ\phi.$$
Applying Lemma \ref{lem_muDP_rewritten}, we see that $u$ is a
solution of the $\muDP$ equation (\ref{muDP}) if and only if
\bea u_t+uu_x&=&-A^{-1}(u(Au)_x-A(uu_x)+3(Au)u_x)\nonumber\\
&=&-A^{-1}(-uu_{xxx}+u_{xx}u_x+uu_{xxx}+2u_xu_{xx}+3(Au)u_x)\nonumber\\
&=&-3A^{-1}(u_xu_{xx}+(Au)u_x)\nonumber\\
&=&-P(u).\nonumber\eea
Recall that
$$\mu(uu_x)=\int_0^1uu_x\dx=\frac{1}{2}\int_0^1\partial_x(u^2)\dx
=\frac{1}{2}(u^2(1)-u^2(0))=0,$$
since $u$ is continuous on $\S^1$. With $u=\xi\circ\phi^{-1}$ the
desired result follows.\endproof
%
%
%
%
\section{Short time existence of geodesics}
We now define the vector field
$$F(\phi,\xi):=(\xi,-P_{\phi}(\xi))$$
such that $(\phi_t,\xi_t)=F(\phi,\xi)$. We know that
$$F:\Diff^n(\S^1)\times\C^n(\S^1)\to\C^n(\S^1)\times\C^n(\S^1),$$
since $P$ is of order zero. We aim to prove smoothness of the map
$F$. It is worth to mention that this will not follow from the
smoothness of $P$ since neither the composition nor the inversion
are smooth maps on $\Diff^n(\S^1)$. The following lemma will be
crucial for our purposes.
\lem\label{lem_polynomial_diffoperator} Assume that $p$ is a
polynomial differential operator of order $r$ with coefficients
depending only on $\mu$, i.e.,
$$p(u)=\sum_{I=(\alpha_0,\ldots,\alpha_r),\atop\alpha_i\in\N\cup\{0\},
\,|I|\leq K}a_I(\mu(u))\,u^{\alpha_0}(u')^{\alpha_1}\cdots
(u^{(r)})^{\alpha_r}.$$
Then the action of $p_{\phi}:=R_{\phi}\circ p\circ R_{\phi^{-1}}$
is
$$p_{\phi}(u)=\sum_{I}a_I\left(\int_{0}^1u(y)\phi_x(y)\d y\right)
q_I(u;\phi_x,\ldots,\phi^{(r)}),$$
where $q_I$ are polynomial differential operators of order $r$
with coefficients being rational functions of the derivatives of $\phi$ up to the order $r$. Moreover, the denominator terms
only depend on $\phi_x$.
\endlem
\proof It is sufficient to consider a monomial
$$m(u)=a(\mu(u))u^{\alpha_0}(u')^{\alpha_1}\cdots(u^{(r)})^{\alpha_r}.$$
We have
$$m_{\phi}(u)=a(\mu(u\circ\phi^{-1}))u^{\alpha_0}[(u\circ\phi^{-1})'
\circ\phi]^{\alpha_1}\cdots[(u\circ\phi^{-1})^{(r)}\circ\phi]^{\alpha_r},$$
where $\circ$ denotes again the composition with respect to the
spatial variable. First, we observe that
$$\mu(u\circ\phi^{-1})=\int_{\S^1}u(\phi^{-1}(x))\d
x=\int_0^1u(y)\phi_x(y)\d y,$$
where we have omitted the time dependence of $u$ and $\phi$.
Recall that $\phi(\S^1)=\S^1$, $\phi_x>0$ and that
$\mu(u\circ\phi^{-1})$ is a constant with respect to the spatial
variable $x\in\S^1$.
Let us introduce the notation
$$a_k=(u\circ\phi^{-1})^{(k)}\circ\phi,\quad k=1,2,\ldots,r.$$
Then, by the chain rule,
$$a_1=(\partial_x(u\circ\phi^{-1}))\circ\phi=
\frac{u_x\circ\phi^{-1}}{\phi_x\circ\phi^{-1}}\circ\phi=\frac{u_x}{\phi_x}$$
and
\bea
a_{k+1}&=&(\partial_x(u\circ\phi^{-1})^{(k)})\circ\phi\nonumber\\
&=&(\partial_x(a_k\circ\phi^{-1}))\circ\phi\nonumber\\
&=&\frac{\partial_xa_k}{\phi_x},\nonumber\eea
so that our theorem follows by induction.\endproof
Recall that in the Banach algebras $\C^n(\S^1)$, $n\geq 1$,
addition and multiplication as well as the mean value operation
$\mu$ and the derivative $\frac{\!\d}{\!\dx}$ are smooth maps. We therefore conclude
that if the coefficients $a_I$ are smooth functions for any
multi-index $I$ and $u$ and $\phi$ are at least $r$ times
continuously differentiable, then $p_{\phi}(u)$ depends smoothly
on $(\phi,u)$.
\prop\label{prop_smooth_vectorfield} The vector field
$$F:\Diff^n(\S^1)\times\C^n(\S^1)\to\C^n(\S^1)\times\C^n(\S^1)$$
is smooth for any $n\geq 3$.
\endprop
\proof We write $F=(F_1,F_2)$. Since $F_1:(\phi,\xi)\mapsto\xi$ is
smooth, it remains to check that $F_2:(\phi,\xi)\mapsto
-P_{\phi}(\xi)$ is smooth. For this purpose, we consider the map
$$\tilde P:\Diff^n(\S^1)\times\C^n(\S^1)\to\Diff^n(\S^1)\times\C^n(\S^1)$$
defined by
$$\tilde P(\phi,\xi)=(\phi,(R_{\phi}\circ P\circ R_{\phi^{-1}})(\xi)).$$
Observe that we have the decomposition $\tilde P=\tilde
A^{-1}\circ\tilde Q$ with
$$\tilde A(\phi,\xi)=(\phi,(R_{\phi}\circ A\circ R_{\phi^{-1}})(\xi))$$
and
$$\tilde Q(\phi,\xi)=(\phi,(R_{\phi}\circ Q\circ R_{\phi^{-1}})(\xi)),$$
where $Q(f):=3(f_xf_{xx}+(Af)f_x)$. We now apply Lemma
\ref{lem_polynomial_diffoperator} to deduce that
$$\tilde A,\tilde Q:\Diff^n(\S^1)\times\C^n(\S^1)\to\Diff^n(\S^1)
\times\C^{n-2}(\S^1)$$
are smooth. To show that $\tilde
A^{-1}:\Diff^n(\S^1)\times\C^{n-2}(\S^1)\to\Diff^n(\S^1)\times
\C^{n}(\S^1)$ is smooth, we compute the derivative $D\tilde A$ at
an arbitrary point $(\phi,\xi)$. We have the following directional
derivatives of the components $\tilde A_1$ and $\tilde A_2$:
$$D_{\phi}\tilde A_1=\id,\quad D_{\xi}\tilde A_1=0,\quad D_{\xi}
\tilde A_2=R_{\phi}\circ A\circ R_{\phi^{-1}}.$$
It remains to compute $(D_{\phi}\tilde
A_2(\phi,\xi))(\psi)=\frac{\!\d}{\!\d\eps}\tilde
A_2(\phi+\eps\psi,\xi)\big|_{\eps=0}$.
In a first step, we calculate
\begin{align}
\partial^2_x(\xi\circ(\phi+\eps\psi)^{-1})
&=\partial_x\left[\left(\frac{\xi_x}{\phi_x+\eps\psi_x}\right)\circ(\phi+\eps\psi)^{-1}\right]\nonumber\\
&=\left(\frac{\xi_{xx}}{(\phi_x+\eps\psi_x)^2}-\xi_{x}\frac{\phi_{xx}+\eps\psi_{xx}}{(\phi_x+\eps\psi_x)^3}\right)\circ(\phi+\eps\psi)^{-1},
\nonumber
\end{align}
from which we get
\begin{align}
\frac{\!\d}{\!\d\eps}\left[
\partial_x^2(\xi\circ(\phi+\eps\psi)^{-1})\circ(\phi+\eps\psi)\right]
&=\frac{\!\d}{\!\d\eps}\left(\frac{\xi_{xx}}{(\phi_x+\eps\psi_x)^2}
-\xi_x\frac{\phi_{xx}+\eps\psi_{xx}}{(\phi_x+\eps\psi_x)^3}\right)\nonumber\\
&=-2\frac{\xi_{xx}\psi_x}{(\phi_x+\eps\psi_x)^3}
-\frac{\xi_x\psi_{xx}}{(\phi_x+\eps\psi_x)^3}\nonumber\\
&\qquad+3\frac{\xi_x\psi_x}{(\phi_x+\eps\psi_x)^4}
(\phi_{xx}+\eps\psi_{xx})\nonumber
\end{align}
and finally
$$\left.\frac{\!\d}{\!\d\eps}\left[\partial_x^2(\xi\circ(\phi+\eps\psi)^{-1})
\circ(\phi+\eps\psi)\right]\right|_{\eps=0}=-2\frac{\xi_{xx}\psi_x}{\phi_x^3}
-\frac{\xi_{x}\psi_{xx}}{\phi_x^3}+3\frac{\phi_{xx}\xi_x\psi_x}{\phi_x^4}.$$
Secondly, we observe that
\bea\left.\frac{\!\d}{\!\d\eps}\mu(\xi\circ(\phi+\eps\psi)^{-1})\right|_{\eps=0}
&=&\left.\frac{\!\d}{\!\d\eps}\int_{\S^1}\xi(y)(\phi_x+\eps\psi_x)(y)
\d y\right|_{\eps=0}\nonumber\\
&=&\int_{\S^1}\xi(y)\psi_x(y)\d y,\nonumber\eea
since $\phi+\eps\psi\in\Diff^n(\S^1)$ for small $\eps>0$.
Hence
$$(D_{\phi}\tilde A_2(\phi,\xi))(\psi)=\int_{\S^1}\xi(y)\psi_x(y)\d y+
2\frac{\xi_{xx}\psi_x}{\phi_x^3}+\frac{\xi_{x}\psi_{xx}}{\phi_x^3}-
3\frac{\phi_{xx}\xi_x\psi_x}{\phi_x^4}$$
and
$$D\tilde A{(\phi,\xi)}=\left(%
\begin{array}{cc}
  \id & 0 \\
  D_{\phi}\tilde A_2(\phi,\xi) & R_{\phi}\circ A\circ R_{\phi^{-1}} \\
\end{array}%
\right).$$
It is easy to check that $D\tilde A{(\phi,\xi)}$ is an invertible
bounded linear operator
$\C^n(\S^1)\times\C^n(\S^1)\to\C^n(\S^1)\times\C^{n-2}(\S^1)$. By
the open mapping theorem, $D\tilde A$ is a topological isomorphism
and, by the inverse mapping theorem, $\tilde A^{-1}$ is
smooth.\endproof
Since $F$ is smooth, we can apply the Banach space version of the
Picard-Lindel\"of Theorem (also known as \emph{Cauchy-Lipschitz
Theorem}) as explained in \cite{Lang}, Chapter XIV-3. This yields
the following theorem about the existence and uniqueness of
integral curves for the vector field $F$.
\thm\label{thm_short_time_ex_of_geo} Given $n\geq 3$,  there is
an open interval $J_n$ centered at zero and an open ball
$B(0,\delta_n)\subset\C^n(\S^1)$ such that for any $u_0\in
B(0,\delta_n)$ there exists a unique solution
$(\phi,\xi)\in\C^{\infty}(J_n,\Diff^n(\S^1)\times\C^n(\S^1))$ of
\mbox{\rm(\ref{muDPgeo})} with initial conditions $\phi(0)=\id$
and $\xi(0)=u_0$. Moreover, the flow $(\phi,\xi)$ depends smoothly
on $(t,u_0)$.\endthm\rm
From Theorem~\ref{thm_short_time_ex_of_geo} we get a unique
short-time solution $u=\xi\circ\phi^{-1}$ of $\muDP$ in
$\C^n(\S^1)$ with continuous dependence on $(t,u_0)$. We now aim
to obtain an analogous result for smooth initial data $u_0$. But since $\C^{\infty}(\S^1)$ is a Fr\'echet space,
classical results like the Picard-Lindel\"of Theorem or the local
inverse theorem for Banach spaces are no longer valid in $\C^{\infty}(\S^1)$. In the
proof of our main theorem, we will make use of a Banach space
approximation of the Fr\'echet space $\C^{\infty}(\S^1)$. First we shall establish that any solution $(\phi,\xi)$ of the $\muDP$ equation
\eqref{muDPgeo} does not lose nor gain spatial regularity as $t$
increases or decreases from zero. For this purpose, the following conservation law is quite useful. In its formulation we use the notation
$m_0(x):=(Au)(0,x)=\mu(u_0)-(u_0)_{xx}$.
%
\lem\label{lem_conservlaw_muDP} Let $u$ be a $\C^3(\S^1)$-solution
of the $\muDP$ equation on $(-T,T)$ and let $\phi$ be the
corresponding flow. Then
$$(Au)(t,\phi(t,x))\phi_x^3(t,x)=m_0,$$
for all $t\in(-T,T)$.
\endlem
\proof We compute
\bea&&\frac{\!\d}{\!\d t}
[(\mu(u)-u_{xx}\circ\phi)\phi_x^3]\nonumber\\
&&=[\mu(u_t)-u_{xxt}\circ\phi-(u_{xxx}\circ\phi)\phi_t]
\,\phi_x^3+3\phi_x^2\phi_{tx}(\mu(u)-u_{xx}\circ\phi)\nonumber\\
&&=[\mu(u_t)-u_{xxt}\circ\phi-(u_{xxx}\circ\phi)(u\circ\phi)]
\,\phi_x^3+3\phi_x^2(u\circ\phi)_x(\mu(u)-u_{xx}\circ\phi)\nonumber\\
&&=[(\mu(u_t)-u_{xxt}-u_{xxx}u)\circ\phi]\,\phi_x^3
+3\phi_x^2(u_x\circ\phi)\phi_x(\mu(u)-u_{xx}\circ\phi)\nonumber\\
&&=[(\mu(u_t)-u_{xxt}-u_{xxx}u)\circ\phi]\,\phi_x^3
+3\phi_x^3[u_x(\mu(u)-u_{xx})]\circ\phi\nonumber\\
&&=[(3u_xu_{xx}-3\mu(u)u_x)\circ\phi]\,\phi_x^3
-3\phi_x^3(u_xu_{xx}-\mu(u)u_x)\circ\phi\nonumber\\
&&=0.\nonumber\eea
Since $\phi(0)=\id$ and $\phi_x(0)=1$, the proof is completed.\endproof
\lem\label{lem_no_gain_or_loss_in_spatial_regularity} Let
$(\phi,\xi)\in\C^{\infty}(J_3,\Diff^3(\S^1)\times\C^3(\S^1))$ be a
solution of $\eqref{muDPgeo}$ with initial data $(\id,u_0)$,
according to Theorem \ref{thm_short_time_ex_of_geo}. Then, for all
$t\in J_3$,
\bea\label{phixx}\phi_{xx}(t)=\phi_x(t)\left(\int_0^t\mu(u)\phi_x(s)\d
s-m_0\int_0^t\phi_x(s)^{-2}\d s\right)\eea
and
\bea\label{xixx}\xi_{xx}(t)
=\xi_x(t)\frac{\phi_{xx}(t)}{\phi_x(t)}+\phi_x(t)
\left[\mu(u)\phi_x(t)-m_0\phi_x(t)^{-2}\right].\eea
\endlem
\proof We have
$$\frac{\!\d}{\!\d t}\left(\frac{\phi_{xx}}{\phi_x}\right)
=\frac{\phi_{xxt}\phi_x-\phi_{xt}\phi_{xx}}{\phi_x^2}.$$
Since $\phi_t=u\circ\varphi$,
$$\phi_{xt}=\phi_{tx}={\partial_x}(u\circ\phi)
=(u_x\circ\phi)\phi_x$$
and
\bea\phi_{xxt}&=&\phi_{txx}\nonumber\\
&=&{\partial_x^2}(u\circ\phi)\nonumber\\
&=&{\partial_x}[(u_x\circ\phi)\phi_x]\nonumber\\
&=&(u_{xx}\circ\phi)\phi_x^2+(u_x\circ\phi)\phi_{xx}.\nonumber\eea
Hence
$$\frac{\!\d}{\!\d t}\left(\frac{\phi_{xx}}{\phi_x}\right)
=(u_{xx}\circ\phi)\phi_x.$$
According to the previous lemma, we know that
$$u_{xx}\circ\phi=\mu(u)-m_0\phi_x^{-3}.$$
Integrating
$$\frac{\!\d}{\!\d t}\left(\frac{\phi_{xx}}{\phi_x}\right)
=\mu(u)\phi_x-m_0\phi_x^{-2}$$
over $[0,t]$ leads to equation (\ref{phixx}) and taking the time
derivative of (\ref{phixx}) yields (\ref{xixx}).\endproof
\rem Since the $\muDP$ equation is equivalent to the quasi-linear
evolution equation
$$u_t+uu_x+3\mu(u)\partial_xA^{-1}u=0,$$
we see that $\mu(u_t)=0$ and hence $\mu(u)=\mu(u_0)$ so that
$\mu(u)$ can in fact be written in front of the first integral
sign in equation (\ref{phixx}).
\endrem
\cor\label{cor_no_gain_or_loss_in_spatial_regularity} Let
$(\phi,\xi)$ be as in Lemma
\ref{lem_no_gain_or_loss_in_spatial_regularity}. If
$u_0\in\C^n(\S^1)$ then we have
$(\phi(t),\xi(t))\in\Diff^n(\S^1)\times\C^n(\S^1)$ for all $t\in
J_3$.
\endcor
\proof We proceed by induction on $n$. For $n=3$ the result is
immediate from our assumption on $(\phi(t),\xi(t))$. Let us assume
that $(\phi(t),\xi(t))\in\Diff^n(\S^1)\times\C^n(\S^1)$ for some
$n\geq 3$. Then Lemma
\ref{lem_no_gain_or_loss_in_spatial_regularity} shows that, if
$u_0\in\C^{n+1}(\S^1)$, then
$(\phi(t),\xi(t))\in\Diff^{n+1}(\S^1)\times\C^{n+1}(\S^1)$,
finishing the proof.
\endproof
\cor\label{cor_no_gain_or_loss_in_spatial_regularity2} Let
$(\phi,\xi)$ be as in Lemma
\ref{lem_no_gain_or_loss_in_spatial_regularity}. If there exists a
nonzero $t\in J_3$ such that $\phi(t)\in\Diff^n(\S^1)$ or
$\xi(t)\in\C^n(\S^1)$ then $\xi(0)=u_0\in\C^n(\S^1)$.
\endcor
\proof Again, we use a recursive argument. For $n=3$, there is
nothing to do. For some $n\geq 3$, suppose that
$u_0\in\C^n(\S^1)$. By the previous corollary,
$(\phi(t),\xi(t))\in\Diff^n(\S^1)\times\C^n(\S^1)$ for all $t\in
J_3$. Assume that there is $0\neq t_0\in J_3$ such that
$\phi(t_0)\in\Diff^{n+1}(\S^1)$ or $\xi(t_0)\in\C^{n+1}(\S^1)$. Since
$\phi_x>0$, Lemma~\ref{lem_no_gain_or_loss_in_spatial_regularity}
immediately implies that also $u_0\in\C^{n+1}(\S^1)$.\endproof
%
Now we discuss Banach space approximations of Fr\'echet spaces.
\defn\label{defn_Banach_space_approx} Let $X$ be a Fr\'echet space.
A \emph{Banach space approximation} of $X$ is a sequence
$\set{(X_n,\norm{\cdot}_n)}{n\in\N_0}$ of Banach spaces such that
$$X_0\supset X_1\supset X_2\supset\cdots\supset X,\quad X=
\bigcap_{n=0}^{\infty}X_n$$
and $\set{\norm{\cdot}_n}{n\in\N_0}$ is a sequence of norms
inducing the topology on $X$ with
$$\norm{x}_0\leq\norm{x}_1\leq\norm{x}_2\leq\ldots$$
for any $x\in X$.\enddefn
We have the following result. For a proof, we refer to
\cite{EscherKolev}.
\lem\label{lem_Banach_space_approx} Let $X$ and $Y$ be Fr\'echet
spaces with Banach space approximations
$\set{(X_n,\norm{\cdot}_n)}{n\in\N_0}$ and
$\set{(Y_n,\norm{\cdot}_n)}{n\in\N_0}$. Let $\Phi_0:U_0\to V_0$ be
a smooth map between the open subsets $U_0\subset X_0$ and
$V_0\subset Y_0$. Let
$$U:=U_0\cap X\quad\text{and}\quad V:=V_0\cap Y,$$
as well as
$$U_n:=U_0\cap X_n\quad\text{and}\quad V_n:=V_0\cap Y_n,$$
for any $n\geq 0$. Furthermore, we assume that, for each $n\geq
0$, the following properties are satisfied:
\begin{enumerate}
\item[\rm (1)] $\Phi_0(U_n)\subset V_n$,%
\item[\rm (2)] the restriction
$\Phi_n:=\left.\Phi_0\right|_{U_n}:U_n\to V_n$ is a smooth map.
\end{enumerate}
Then $\Phi_0(U)\subset V$ and the map
$\Phi:=\left.\Phi_0\right|_{U}:U\to V$ is smooth.\endlem\rm
Now we come to our main theorem which we first formulate in the
geometric picture.
\thm\label{thm_main_Escher_Kolev09} There exists an open interval
$J$ centered at zero and $\delta>0$ such that for all
$u_0\in\C^{\infty}(\S^1)$ with $\norm{u_0}_{\mbox{\rm\scriptsize
C}^3(\S^1)}<\delta$, there exists a unique solution
$(\phi,\xi)\in\C^{\infty}(J,\Diff^{\infty}(\S^1)\times\C^{\infty}(\S^1))$
of \mbox{\rm(\ref{muDPgeo})} such that $\phi(0)=\id$ and
$\xi(0)=u_0$. Moreover, the flow $(\phi,\xi)$ depends smoothly on
$(t,u_0)\in J\times\C^{\infty}(\S^1)$.
\endthm
\proof Theorem \ref{thm_short_time_ex_of_geo} for $n=3$ shows that
there is an open interval $J$ centered at zero and an open ball
$U_3=B(0,\delta)\subset\C^3(\S^1)$ so that for any $u_0\in U_3$
there exists a unique solution
$(\phi,\xi)\in\C^{\infty}(J,\Diff^{3}(\S^1)\times\C^{3}(\S^1))$ of
(\ref{muDPgeo}) with initial data $(\id,u_0)$ and a smooth flow
$$\Phi_3:J\times U_3\to\Diff^3(\S^1)\times\C^3(\S^1).$$
Let
$$U_n:=U_3\cap\C^n(\S^1)\quad\text{and}\quad U_{\infty}
:=U_3\cap\C^{\infty}(\S^1).$$
By Corollary \ref{cor_no_gain_or_loss_in_spatial_regularity}, we
have
$$\Phi_3(J\times U_n)\subset\Diff^n(\S^1)\times\C^n(\S^1)$$
for any $n\geq 3$ and the map
$$\Phi_n:=\left.\Phi_3\right|_{J\times U_n}:
J\times U_n\to\Diff^n(\S^1)\times\C^n(\S^1)$$
is smooth. Lemma \ref{lem_Banach_space_approx} implies that
$$\Phi_3(J\times U_{\infty})\subset
\Diff^{\infty}(\S^1)\times\C^{\infty}(\S^1),$$
completing the proof of the short-time existence for smooth initial data $u_0$.
Moreover, the mapping
$$\Phi_{\infty}:=\left.\Phi_3\right|_{J\times U_{\infty}}:
J\times U_{\infty}\to\Diff^{\infty}(\S^1)\times\C^{\infty}(\S^1)$$
is smooth, proving the smooth dependence on time and on the initial
condition.\endproof
Under the assumptions of Theorem \ref{thm_main_Escher_Kolev09},
the map
$$\Diff^{\infty}(\S^1)\times\C^{\infty}(\S^1)\to
\C^{\infty}(\S^1),\quad(\phi,\xi)\mapsto \xi\circ\phi^{-1}=u$$
is smooth. Thus we obtain the result stated in Theorem
\ref{cor_main_Escher_Kolev_09}.
\section{The exponential map}
For a Banach manifold $M$ equipped with a symmetric linear
connection, the exponential map is defined as the time one of the
geodesic flow, i.e., if $t\mapsto\gamma(t)$ is the (unique)
geodesic in $M$ starting at $p=\gamma(0)$ with velocity
$\gamma_t(0)=u\in T_pM$ then $\exp_p(u)=\gamma(1)$. Roughly
speaking, the map $\exp_p(\cdot)$ is a projection from $T_pM$ to
the manifold $M$. Since the derivative of $\exp_p$ at zero is the
identity, the exponential map is a smooth diffeomorphism from a
neighbourhood of zero of $T_pM$ to a neighbourhood of $p\in M$.
However, this fails for Fr\'echet manifolds like
$\Diff^{\infty}(\S^1)$ in general. We know that the Riemannian
exponential map for the $L^2$-metric on $\Diff^{\infty}(\S^1)$ is
not a local $\C^1$-diffeomorphism near the origin, cf.\
\cite{ConstantinKolev2}. For the Camassa-Holm equation and more
general for the $H^k$-metrics, $k\geq 1$, the Riemannian
exponential map in fact is a smooth local diffeomorphism. This
result was generalized to the family of $b$-equations, see
\cite{EscherKolev}, and in this section we obtain a similar result
for the $\muDP$ equation.\\\indent
The basic idea of the proof of Theorem \ref{thm_exp} is to
consider a perturbed problem: Let $(\phi^{\eps},\xi^{\eps})$
denote the local expression of an integral curve of
(\ref{muDPgeo}) in $T\Diff^n(\S^1)$ with initial data $(\id,u+\eps
w)$, where $u,w\in\C^n(\S^1)$. Let
$$\psi(t):=\left.\frac{\partial\phi^{\eps}(t)}
{\partial\eps}\right|_{\eps=0}.$$
By the homogeneity of the geodesics,
$$\phi^{\eps}(t)=\exp(t(u+\eps w)),$$
so that
$$\psi(t)=D\left(\exp(tu)\right)tw=:L_n(t,u)w,$$
where $L_n(t,u)$ is a bounded linear operator on $C^n(\S^1)$.
\lem Suppose that $u\in\C^{n+1}(\S^1)$. Then, for $t\neq 0$,
$$L_n(t,u)(\C^n(\S^1)\backslash\C^{n+1}(\S^1))\subset
\C^n(\S^1)\backslash\C^{n+1}(\S^1).$$
\proof First, we write down equation (\ref{phixx}) for
$\phi^{\eps}(t)$,
$$\phi_{xx}^{\eps}(t)=\phi_{x}^{\eps}(t)\left[\mu(u+\eps w)\int_0^t
\phi_x^{\eps}(s)\d
s-m_0^{\eps}\int_0^t\phi_x^{\eps}(s)^{-2}\d s\right],$$
and take the derivative with respect to $\eps$,
\begin{align}
\frac{\partial\phi_{xx}^{\eps}}{\partial\eps}(t)
&=\frac{\partial\phi_{x}^{\eps}}{\partial\eps}(t) \left[\mu(u+\eps
w)\int_0^t\phi_x^{\eps}(s)\d s
-m_0^{\eps}\int_0^t\phi_x^{\eps}(s)^{-2}\d s\right]\nonumber\\
&\qquad+\phi_x^{\eps}(t)\left[\mu(w)\int_0^t\phi_x^{\eps}(s)\d
s+\mu(u+\eps
w)\int_0^t\frac{\partial\phi_x^{\eps}}{\partial\eps}(s)\d
s\right]\nonumber\\
&\qquad-\phi_x^{\eps}(t)\left[\frac{\partial
m_0^{\eps}}{\partial\eps} \int_0^t\phi_x^{\eps}(s)^{-2}\d
s+m_0^{\eps}\int_0^t
\frac{\partial}{\partial\eps}\phi_x^{\eps}(s)^{-2}\d
s\right].\nonumber
\end{align}
Notice that
$$\frac{\partial m_0^{\eps}}{\partial\eps}=\mu(w)-w_{xx}=Aw$$
and that $m_0^{\eps}\to m_0=Au$ as $\eps\to 0$.
Hence
\begin{align}
\psi_{xx}(t)&=\psi_x(t)\left[\mu(u)\int_0^t\phi_x(s)\d
s-m_0\int_0^t\phi_x(s)^{-2}\d s\right]\nonumber\\
&\quad+\phi_x(t)\left[\mu(w)\int_0^t\phi_x(s)\d
s+\mu(u)\int_0^t\psi_x(s)\d s\right]\nonumber\\
&\quad-\phi_x(t)\bigg[(\mu(w)-w_{xx})\int_0^t\phi_x(x)^{-2}\d
s-2m_0\int_0^t\psi_x(s)\phi_x(s)^{-3}\d
s\bigg]\nonumber\\
&=a(t)\psi_x(t)+b(t)\int_0^tc(s)\psi_x(s)\d
s+d(t)+e(t)w_{xx}\nonumber
\end{align}
with $a(t),b(t),c(t),d(t),e(t)\in\C^{n-1}(\S^1)$ and $e(t)\neq 0$
for $t\neq 0$. Finally, if
$$w\in\C^{n}(\S^1)\backslash\C^{n+1}(\S^1),$$
then
$$\psi(t)=L_n(t,u)w\in\C^{n}(\S^1)\backslash\C^{n+1}(\S^1).$$
\endproof
Let us now turn to the proof of Theorem \ref{thm_exp}. Since
$\C^3(\S^1)$ is a Banach space and $\Diff^3(\S^1)$ is a Banach
manifold modelled over $\C^3(\S^1)$, we know that the exponential
map is a smooth diffeomorphism near zero, i.e., there are
neighbourhoods $U_3$ of zero in $\C^3(\S^1)$ and $V_3$ of $\id$ in
$\Diff^3(\S^1)$ such that
$$\exp_3:=\exp|_{U_3}:U_3\to V_3$$
is a smooth diffeomorphism. For $n\geq 3$, we now define
$$U_n:=U_3\cap\C^n(\S^1)\quad\text{and}\quad V_n=V_3\cap\Diff^n(\S^1).$$
Let $\exp_n:=\exp_3|_{U_n}$. Since $\exp_n$ is a restriction of
$\exp_3$, it is clearly injective. We now use Corollary
\ref{cor_no_gain_or_loss_in_spatial_regularity} and Corollary
\ref{cor_no_gain_or_loss_in_spatial_regularity2} to deduce that
$\exp_n$ is also surjective, more precisely, $\exp_n(U_n)=V_n$. If
the geodesic $\phi$ with $\phi(1)=\exp(u)$ starts at
$\id\in\Diff^n(\S^1)$ with velocity vector $u$ belonging to
$\C^{n}(\S^1)$, then $\phi(t)\in\Diff^n(\S^1)$ for any $t$ and
hence $\exp_n(U_n)\subset V_n$. Conversely, if $v\in V_n$ is
given, then there is $u\in U_3$ with $\exp_3(u)=v$. Corollary
\ref{cor_no_gain_or_loss_in_spatial_regularity2} immediately
implies that $u\in\C^n(\S^1)$; hence $u\in U_n$ and $\exp_n(u)=v$.
Note that $\exp_n$ is a bijection from $U_n$ to $V_n$.
Furthermore, $\exp_n$ is a smooth map and diffeomorphic $U_n\to V_n$. We now show
that $\exp_n$ is a smooth diffeomorphism; precisely we show that
$\exp_n^{-1}:V_n\to U_n$ is smooth by virtue of the inverse
mapping theorem. For each $u\in\C^n(\S^1)$, $D\exp_n(u)$ is a
bounded linear operator $\C^n(\S^1)\to\C^n(\S^1)$. Notice that
$$D\exp_n(u)=D\exp_3(u)|_{\text{C}^n(\S^1)},$$
from which we conclude that $D\exp_n(u)$ is injective. Let us
prove the surjectivity of $D\exp_n(u)$, $n\geq 3$, by induction.
For $n=3$, this follows from the fact that $\exp_3:U_3\to V_3$ is
diffeomorphic and hence a submersion. Assume that $D\exp_n(u)$ is
surjective for some $n\geq 3$ and that $u\in\C^{n+1}(\S^1)$. We
have to show that this implies the surjectivity of
$D\exp_{n+1}(u)$. But this is a direct consequence of
$D\exp_n(u)=L_n(1,u)$ and the previous lemma: Let $f\in
\C^{n+1}(\S^1)$. We have to find
$g\in\C^{n+1}(\S^1)$ with the property $D\exp_{n+1}(u)g=f$. By our
assumption, there is $g\in\C^n(\S^1)$ such that $D\exp_n(u)g=f$.
Suppose that $g\notin\C^{n+1}(\S^1)$. But then
$f=L_n(1,u)g\notin\C^{n+1}(\S^1)$ in contradiction to the choice
of $f$. Thus $g\in\C^{n+1}(\S^1)$ and $D\exp_{n+1}(u)g=f$. Now we
can apply the open mapping theorem to deduce that for any $n\geq
3$ and any $u\in\C^n(\S^1)$ the map
$$D\exp_n(u):\C^n(\S^1)\to\C^n(\S^1)$$
is a topological isomorphism. By the inverse function theorem,
$\exp_n:U_n\to V_n$ is a smooth diffeomorphism. If we define
$$U_{\infty}:=U_3\cap\C^{\infty}(\S^1)\quad\text{and}
\quad V_{\infty}:=V_3\cap\Diff^{\infty}(\S^1),$$
Lemma \ref{lem_Banach_space_approx} yields that
$$\exp_{\infty}:=\exp_3|_{U_{\infty}}:U_{\infty}\to V_{\infty}$$
as well as
$$\exp_{\infty}^{-1}:V_{\infty}\to U_{\infty}$$
are smooth maps. Thus $\exp_{\infty}$ is a smooth diffeomorphism
between $U_{\infty}$ and $V_{\infty}$.
%
%
%
%
\setlength{\bibsep}{0.1\baselineskip}
\small\bibliography{Escher_Kohlmann_Kolev_final_lit}
\end{document}